\documentclass[10.5pt]{article}
\usepackage{amsfonts,mathrsfs}

\def\sqr#1#2{{\vcenter{\vbox{\hrule height.#2pt
              \hbox{\vrule width.#2pt height#1pt \kern#1pt \vrule width.#2pt}
              \hrule height.#2pt}}}}
\def\3n{\negthinspace \negthinspace \negthinspace }
\def\2n{\negthinspace \negthinspace }
\def\1n{\negthinspace }

\def\dbE{\mathbb{E}}
\def\dbF{\mathbb{F}}

\def\dbP{\mathbb{P}}

\def\dbR{\mathbb{R}}
\def\dbS{\mathbb{S}}

\def\sD{\mathscr{D}}

\def\sG{\mathscr{G}}
\def\sH{\mathscr{H}}

\def\sL{\mathscr{L}}
\def\sM{\mathscr{M}}
\def\sN{\mathscr{N}}

\def\sP{\mathscr{P}}
\def\sQ{\mathscr{Q}}
\def\sR{\mathscr{R}}
\def\sS{\mathscr{S}}

\def\sU{\mathscr{U}}

\def\sX{\mathscr{X}}

\def\={\buildrel \triangle \over =}

\def\ds{\displaystyle}

\def\ns{\noalign{\ss}}

\def\e{\varepsilon}

\def\si{\sigma}

\def\f{\varphi}
\def\th{\theta}

\def\D{\Delta}
\def\Th{\Theta}
\def\L{\Lambda}

\def\cF{{\cal F}}

\def\cL{{\cal L}}
\def\cM{{\cal M}}
\def\cN{{\cal N}}

\def\cP{{\cal P}}

\def\cR{{\cal R}}

\def\cW{{\cal W}}

\def\mds{\medskip}

\def\ss{\smallskip}
\def\ms{\medskip}

\def\q{\quad}
\def\qq{\qquad}
\def\hb{\hbox}

\def\lan{\mathop{\langle}}
\def\ran{\mathop{\rangle}}

\def\h{\widehat}
\def\wt{\widetilde}
\def\cd{\cdot}

\def\deq{\mathop{\buildrel\D\over=}}

\def\({\Big (}
\def\){\Big )}
\def\[{\Big[}
\def\]{\Big]}
\def\bde{\begin{definition}}
\def\ede{\end{definition}}
\def\be{\begin{equation}}
\def\bel{\begin{equation}\label}
\def\ee{\end{equation}}
\def\bt{\begin{theorem}}
\def\et{\end{theorem}}
\def\bc{\begin{corollary}}
\def\ec{\end{corollary}}
\def\bl{\begin{lemma}}
\def\el{\end{lemma}}
\def\bp{\begin{proposition}}
\def\ep{\end{proposition}}
\def\bas{\begin{assumption}}
\def\eas{\end{assumption}}
\def\br{\begin{remark}}
\def\er{\end{remark}}
\def\ba{\begin{array}}
\def\ea{\end{array}}
\def\ed{\end{document}}

\def\square#1{\vbox{\hrule\hbox{\vrule height#1%
     \kern#1\vrule}\hrule}}
\def\rectangle#1#2{\vbox{\hrule\hbox{\vrule height#1%
     \kern#2\vrule}\hrule}}

\font\tenbb=msbm10 \font\sevenbb=msbm7 \font\fivebb=msbm5

\newfam\bbfam
\scriptscriptfont\bbfam=\fivebb \textfont\bbfam=\tenbb
\scriptfont\bbfam=\sevenbb

\newtheorem{lemma}{Lemma}[section]
\newtheorem{remark}{Remark}[section]

\newtheorem{theorem}{Theorem}[section]
\newtheorem{corollary}{Corollary}[section]

\newtheorem{definition}{Definition}[section]
\newtheorem{proposition}{Proposition}[section]
\newtheorem{assumption}{Assumption}[section]
\newtheorem{proof}{Proof}[section]

\makeatletter
   
   \@addtoreset{equation}{section}
\makeatother

\usepackage{amsfonts}
\usepackage{latexsym}
\usepackage{amsmath}
\usepackage{amssymb}
\usepackage{color}

\begin{document}

\title{Equilibrium controls in time inconsistent stochastic linear quadratic
problems\footnote{The research was supported by the NSF of China under grant 11231007, 11401404 and 11471231.}}

\author{Tianxiao Wang \footnote{School of
Mathematics, Sichuan University, Chengdu, P. R. China. Email:wtxiao2014@scu.edu.cn.}}

\maketitle

\begin{abstract}
This paper deals with a class of time inconsistent stochastic linear quadratic (SLQ) optimal control problems in Markovian framework. Three notions, i.e., closed-loop equilibrium controls/strategies, open-loop equilibrium controls and their closed-loop representations, are characterized in $unified$ manners. These results indicate clearer and deeper distinctions among these notions. For example, in particular time consistent setting, the open-loop equilibrium controls are fully characterized by $first$-$order$, $second$-$order$ $necessary$ $optimality$ $conditions$, and become needlessly optimal, while the closed-loop equilibrium controls naturally reduce into $closed$-$loop$ $optimal$ $controls$.

\end{abstract}

\ms


\ms



\section{Introduction}
Through out this paper, $(\Omega,\cF,\dbP,\dbF)$ is a complete filtered probability space, on which one-dimensional standard Brownian motion $W(\cd)$ is
defined. Here $\dbF\equiv\{\cF_t\}_{t\ge0}$ is the natural filtration of $W(\cd)$ augmented by $\dbP$-null sets.

\subsection{Formulation of time inconsistent optimal control problems}

For any $t\in[0,T)$, we consider the
following stochastic differential equation (SDE):
\bel{state-equation-original}\left\{\2n\ba{ll}
\ns\ds
dX(s)=\big[A(s)X(s)+B(s)u(s)+ b(s) \big]ds\\
\ns\ds\qq\qq +\big[C(s)X(s)+D(s)u(s)+ \si(s) \big]dW(s),\q
s\in[t,T],\\
\ns\ds X(t)=\xi,\ea\right.\ee
and the cost functional defined by
\bel{cost-original-classical-static}\ba{ll}
\ns\ds J(t,\xi;u(\cd))={1\over2}\dbE_t\Big\{\int_t^T\big[\lan
Q(s)X(s),X(s)\ran+2\lan S(s)X(s),u(s)\ran\\
\ns\ds\qq\qq\qq+\lan
R(s)u(s),u(s)\ran\big]ds +\lan GX(T),X(T)\ran\Big\}.\ea\ee
Here $A,B,C,D,Q, S, R,G$ are suitable matrix-valued
(deterministic) functions,
$b,\si$ are proper stochastic processes, and $\dbE_t(\cd):=\dbE[\,\cd\,|\cF_t]$ stands for conditional expectation operator. In the above, $X(\cd)$, valued
in $\dbR^n$, is called the {\it state process}, $u(\cd)$, valued in
$\dbR^m$, is called the {\it control process}, and $(t,\xi)\in\sD$ is
called the {\it initial pair} where
$$\sD:=\Big\{(t,\xi)\bigm|t\in[0,T],~\xi\hb{ is $\cF_t$-measurable, }\dbE|\xi|^2<\infty\Big\}.$$
We denote the set of all control processes by
$$\ba{ll}
\ds\sU[t,T]\equiv \Big\{u:[t,T]\times\Omega\to\dbR^m\bigm|u
\hb{ is $\dbF$-progressively measurable},\\
\ns\ds\qq\qq\qq\qq\qq\qq\qq\ \ \dbE\int_t^T|u(s)|^2ds<\infty\Big\}.
\ea$$
Under some mild
conditions on the coefficients, for any initial pair $(t,\xi)$
and a control $u(\cd)\in\sU[t,T]$, the state equation (\ref{state-equation-original})
admits a unique solution $X(\cd)=X(\cd\,;t,x,u(\cd))$, and the cost
functional $J(t,\xi;u(\cd))$ is well-defined. We pose the
following stochastic linear quadratic (SLQ) optimal control problem.

\medskip

\bf Problem (SLQ). \rm For any given $(t,\xi)$, find
a $\bar u(\cd)\in\sU[t,T]$ such that
\bel{1.3}J(t,\xi;\bar
u(\cd))=\inf_{u(\cd)\in\sU[t,T]}J(t,\xi;u(\cd))\deq V(t,\xi).\ee
Any $\bar u(\cd)\in\sU[t,T]$ satisfying (\ref{1.3}) is called an
{\it optimal control} for the given initial pair $(t,\xi)$, the
corresponding state process $\bar X(\cd)$ is called an {\it optimal
state process} for $(t,\xi)$, $(\bar X(\cd),\bar u(\cd))$ is called an
{\it optimal pair} for $(t,\xi)$, and $V(\cd\,,\cd)$ is called the
{\it value function} of Problem (SLQ).

\medskip

For above optimal control problem, it is reasonable to keep the state process stable with respect to possible variation of random factors. To this end, one effective way is to add the variation of $X(\cd)$, i.e.
$$\hbox{Var}_t[X]:=\dbE_t\big[X(T)-\dbE_t X(T)\big]^2=\dbE_t |X(T)|^2-\big[\dbE_t X(T)\big]^2
$$
into the cost functional (e.g., \cite{Bjork-Murgoci-2014}, \cite{Bjork-Murgoci-Zhou-2012}, \cite{Hu-Jin-Zhou-2012}, \cite{Hu-Jin-Zhou-2017}, \cite{Huang-Li-Wang-2017}, \cite{Li-Li-2013}, \cite{Wei-Wang-2017}, \cite{Yong-2017}, etc). Therefore, it is natural to propose the following general modified cost functional
\bel{cost-original}\ba{ll}
\ns\ds J(t,\xi;u(\cd)) ={1\over2}\dbE_t\Big\{\int_t^T\big[\lan
Q(s)X(s),X(s)\ran+2\lan S(s)X(s),u(s)\ran\\
\ns\ds\qq\qq\qq+\lan\wt
Q(s)\dbE_t[X(s)],\dbE_t[X(s)]\ran+2\lan\wt
S(s)\dbE_t[X(s)],\dbE_t[u(s)]\ran\\
\ns\ds\qq\qq\qq +\lan
R(s)u(s),u(s)\ran +\lan\wt
R(s)\dbE_t[u(s)],\dbE_t[u(s)]\ran\big]ds\\
\ns\ds\qq\qq\qq+\lan GX(T),X(T)\ran+\lan\wt
G\dbE_t[X(T)],\dbE_t[X(T)]\ran+2\lan g,\dbE_tX(T)\ran\Big\}.\ea\ee
Here $\wt S,\wt R, \wt G,\wt Q$ are deterministic matrices-valued functions and $g$ is a vector.

In this scenario, the optimal controls become time-inconsistent, i.e.,
the ``optimal'' control based on this moment may not
keep optimality in future. We refer to \cite{Yong-2017} for some explicit examples.

\subsection{Related literature}

The study on time inconsistency by economists
actually dates back to Strotz [12] in the 1950s. One
possible way to treat time inconsistency is to discuss the pre-committed controls for
which the solutions are verified to be optimal only at
the initial time.

In this paper, we shall discuss above optimal control problem from
another viewpoint. More precisely, we investigate the time inconsistency
within a game-theoretic framework and analyze
the time-consistent equilibrium solution (e.g., \cite{Pollak-1968}, \cite{Peleg-Yaari-1973}, \cite{Goldman-1980}). Recently, people began to treat the equilibrium controls using the ideas of stochastic control theories, and developed several different approaches in the existing papers. These methods range from dynamic programming principles and verification procedures to maximum principles and variational techniques.

$\diamond$ In Bj\"{o}rk-Murgoci \cite{Bjork-Murgoci-2012}, Bj\"{o}rk et al \cite{Bjork-Murgoci-2014}, the
authors examined a general class of time inconsistent problems under
Markovian framework by equilibrium
value functions. In the continuous case, they formally derived the extended HJB equations, and then rigorously proved the verification theorem by the conclusions of discrete time case, see Theorem 5.2 in \cite{Bjork-Murgoci-2014}. They also present some special cases including a linear
quadratic control problem in which equilibrium solutions are constructed. This method was also used to treat investment-reinsurance problems with mean-variance criterion, see e.g., \cite{Li-Li-2013}, \cite{Zeng-Li-2011}.

$\diamond$  In Yong (\cite{Yong-2012-2}, \cite{Yong-2017}), the author discussed a class of time inconsistent optimal control problems by multi-person differential games approach, where a new kind of equilibrium HJB equations/sytems of Riccati equations were introduced. Unlike \cite{Bjork-Murgoci-2012}, \cite{Bjork-Murgoci-2014}, they started the investigations in continuous time setting, made partition on time intervals and used tricks of
forward-backward stochastic differential equations (FBSDEs). Further study along this can be found in \cite{Wang-Wu-2016}, \cite{Wei-Yong-Yu-2017}, and so on.

$\diamond$ In Ekeland and Lazrak (\cite{Ekeland-Mbodji-Pirvu-2012}, \cite{Ekeland-Pirvu-2008}), they considered some financial problems such as investment and consumption model with time-inconsistency feature. They used the variational ideas to introduce certain feedback/closed-loop equilibrium controls, and spread out discussions via equilibrium value functions. Compared with the general situation in \cite{Bjork-Murgoci-2012}, \cite{Bjork-Murgoci-2014}, the particular form of equilibrium value functions were proposed according to the given cost functional, while the complex convergence arguments were avoided.

$\diamond$ Inspired by the ideas of stochastic maximum principles in optimal control theories, Hu et al. \cite{Hu-Jin-Zhou-2012} studied a class of time inconsistent SLQ problems in Markovian setting, introduced
open-loop equilibrium controls and their closed-loop representations, derived general sufficient conditions through a flow of FBSDEs or systems of backward ordinary differential equations (ODEs). Just recently, the same authors continued to discuss the uniqueness of open-loop equilibrium controls in \cite{Hu-Jin-Zhou-2017}. More related details can also be found in \cite{Djehiche-Huang-2016}, \cite{Wei-Wang-2017}, \cite{Wang-Wu-2015}.

\subsection{Unified approach and contributions}

As to Problem (SLQ), in this article we propose a unified method to characterize the open-loop
equilibrium controls, the closed-loop representations of open-loop equilibrium controls, closed-loop equilibrium controls/strategies. We combines the ideas from variational analysis, forward-backward stochastic differential equations and forward-backward decoupling procedures. In the following, we provide a brief outline of our approach.

For any $(\Th_1,\Th_2,\f)\in L^2(0,T;\dbR^{m\times n})\times L^2(0,T;\dbR^{m\times n})\times L^2_{\dbF}(0,T;\dbR^m),$ we start with control processes
\bel{Definition-perturbed-introduction}\ba{ll}
\ns\ds
u:=(\Th_1+\Th_2)X+\f,\ \ u^\e:=\Th_1 X^\e+\Th_2 X+\f+vI_{[t,t+\e]}.
\ea\ee
They can reduce into the required equilibrium controls and perturbed controls in various settings. More precisely, if $\Th_2\equiv0$, or $\Th_1\equiv0$, or $\Th_1\equiv\Th_2\equiv0$, $u$ and $u^\e$ play the important roles in obtaining closed-loop equilibrium controls/strategies, open-loop equilibrium controls, the closed-loop representation of open-loop equilibrium controls, respectively. We refer to Subsection \ref{Main-proof-subsection} for more detailed discussions.

In view of the definitions for equilibrium controls, we proceed to consider the difference of the cost functional at $u$, $u^\e$. To do so, given $X$ and $X^\e$,
we introduce, respectively, backward stochastic differential equations (BSDEs) with conditional expectations. We point out that the one associated with $X^\e$ appears for the first time in the literature. As a result, we obtain two forward-backward systems in which the terminal parts and generators of backward systems rely respectively on $X$, $X^\e$.

To tackle the limit part in the definitions of both open-loop and closed-loop equilibrium controls (i.e., Definitions \ref{Definition-1}, Definition \ref{Definition-3} next), we continue to decouple the above two forward-backward systems. More precisely, we make conjectures on the solutions of backward systems, formally obtain a class of systems of BSDEs merely depending on given coefficients, and then verify our arguments rigorously. At last we establish our characterizations with proper convergence procedures.

At this very moment, it is worth mentioning that the previous proposed approach demonstrates several new advantages on the treatment of both open-loop equilibrium controls, closed-loop equilibrium controls/strategies. Unlike \cite{Bjork-Murgoci-2012}, \cite{Bjork-Murgoci-2014}, \cite{Yong-2012-2}, \cite{Yong-2017}, our procedures on closed-loop equilibrium strategy in continuous time drop the reliance on complex convergence arguments from discrete time to continuous case. Comparing with \cite{Hu-Jin-Zhou-2012}, \cite{Hu-Jin-Zhou-2017}, our methodology on open-loop equilibrium controls neither requires any non-definite assumptions on the involved coefficients, nor directly uses the conclusions of stochastic maximum principles. Moreover, it can be adjusted into the random coefficients case, see \cite{Wei-Wang-2017}.

Even though both open-loop equilibrium controls and closed-loop equilibrium controls are widely investigated in the literature, there is no paper discussing their differences to our best. In this paper, we give a clear picture by the obtained characterizations. For example, in the classical SLQ setting, open-loop equilibrium controls are fully characterized by first-order, second-order necessary conditions. In other words, they are weaker than optimal controls ( Remark \ref{Relation-between-open-loop-closed-loop-2}). However, in the same situation, the closed-loop equilibrium controls happen to reduce exactly into closed-loop optimal controls (Remark \ref{Relation-between-open-loop-closed-loop-2}). Eventually, we point out that the characterizations on open-loop, closed-loop equilibrium controls, respectively, include two different $second$-$order$ $equilibrium$ $conditions$, which are absent in nearly all the relevant articles.

%
%

\subsection{Outline of the article}

The remainder of this article of structured as follows. In Section 2, an overview of assumptions, notation used in the sequel is provided. In Section 3, the main conclusions of this article are gathered and some important remarks are demonstrated. In Section 4, the proofs of the main results in Section 3 are given. Section 5 concludes this article.

\section{Preliminary notations}

Given $H:=\dbR^n,\dbR^{n\times n},\dbS^{n\times n},$ etc, we introduce the following hypotheses on the coefficients of (\ref{state-equation-original}), (\ref{cost-original}).

\medskip

\bf (H1) \rm Suppose $A, \ B,\ C,\ D, \ R,\ \wt R, \ Q,\ \wt Q,\ S,\ \wt S \in L^{\infty}(0,T;H),$    $G,\ \wt G,\ g\in H,$ $b \in L^2_{\dbF}(\Omega;L^1(0,T;H)),$ $\si\in L^2_{\dbF}(0,T;H)$.

\medskip

For $0\leq s\leq t\leq T$, we also define some involved spaces as follows.
$$\ba{ll}
\ns\ds L^2_{\dbF}(s,t;H):= \Big\{X:[s,t]\times\Omega\to
H\Bigm|X(\cdot) \hbox{ is
$\dbF$-adapted, measurable,\ }\\
\ns\ds\qq\qq\qq\qq\qq\qq\qq\qq \dbE\int_s^t|X(r)|^2dr<\infty\Big\},\\
\ns\ds L^{\infty}(s,t; H):= \Big\{X:[s,t]\to H\Bigm|X \hbox{ is
deterministic, measurable,}\ \sup\limits_{r\in
[s,t]}|X(r)|<\infty\Big\},\\
\ns\ds L^2_{\dbF}(\Omega;L^1(s,t;H)):= \Big\{X : [s,t]\times\Omega\rightarrow
H\bigm|X(\cd)\hbox{ is $\dbF$-adapted, measurable,}\\
\ns\ds\qq\qq\qq\qq\qq\qq\qq\qq\qq\q  \dbE\Big[\ds\int_s^t|X(r)|dr\Big]^2<\infty
 \Big\},\\
\ns\ds L^2_{\dbF}(\Omega;C([s,t];H)):= \Big\{X : [s,t]\times\Omega\rightarrow
H\bigm|X(\cd)\hbox{ is $\dbF$-adapted, measurable}\\
\ns\ds \qq\qq\qq\qq\qq\qq\qq\qq\q \qq\hbox{continuous}\ \dbE\sup\limits_{r\in[s,t]}|X(r)|^2<\infty\Big\}.
\ea$$
To begin with, we look at Problem (SLQ) from an
open-loop equilibrium control viewpoint. The following definition is adapted from \cite{Hu-Jin-Zhou-2012}, \cite{Hu-Jin-Zhou-2017}.

\medskip

\begin{definition}\label{Definition-1}
Given $X^*(0)=x_0\in\dbR^n$, a state-control pair $(X^*,u^*)\in L^2_{\dbF}(\Omega;C([0,T];\dbR^n))\times L^2_{\dbF}(0,T;\dbR^m)$
is called an {\it open-loop equilibrium pair} if for any
$t\in[0,T)$, small $\e>0$, $\mathcal{F}_t$-measurable $v$ satisfying
$\dbE|v|^2<\infty$, the following holds:
\bel{optimal-open}\lim_{\overline{\e\to0}}
{J(t,X^*(t);u^{v,\e}(\cd))-J\big(t,X^*(t);u^*(\cd)\big|_{[t,T]}\big)
\over\e}\ge0,\ee
where $u^{v,\e} =u^*+vI_{[t,t+\e]}$. Here $u^*$ and $X^*$ are called {\it open-loop equilibrium
control} and {\it open-loop equilibrium state process}.
\end{definition}

\medskip

Roughly speaking, the definition shows the $dynamic$ $local$ $optimality$ in some manner. In this paper we will explore deeper properties of such equilibrium controls via their characterizations.

Due to our particular linear quadratic structure, we also introduce the closed-loop representation of open-loop equilibrium control $u^*$.

\medskip

\begin{definition}\label{Definition-2}
An open-loop equilibrium control $u^*\in L^2_{\dbF}(0,T;\dbR^m)$ associated with $X^*(0)=x_0\in\dbR^n$ is said to have a $closed$-$loop$ $representation$ if $u^*=\Th^*X^*+\f^*$ where $X^*$ is the associated state process on $[0,T]$, and $(\Th^*,\f^*)\in L^2(0,T;\dbR^{m\times n})\times L^2_{\dbF}(0,T;\dbR^{m})$. Here they are called $open$-$loop$ $equilibrium$ $strategy$ pair, which are independent of $x_0$.
\end{definition}

\medskip

From the open-loop strategy viewpoint, we can capture more explicit expression of open-loop equilibrium control. However, this kind of strategy is distinctive from the following one.

\medskip

\begin{definition}\label{Definition-3}
$(\Th^*,\f^*)\in L^2(0,T;\dbR^{m\times m}) \times L^2 _{\dbF}(0,T;\dbR^{m}) $ is called a  {\it closed-loop equilibrium strategy}, if for any initial state $x_0\in\dbR^n$,
$t\in[0,T)$, small $\e>0$, $\mathcal{F}_t$-measurable $v$ satisfying
$\dbE|v|^2<\infty$,
\bel{optimal-closed}\lim_{\overline{\e\to0}}
{J(t,X^*(t);u^{\e}(\cd))-J\big(t,X^*(t);u^*(\cd)\big|_{[t,T]}\big)
\over\e}\ge0,\ee
where
$u^*:=\Th^*X +\f^*,$ $u^{\e}:=\Th^*X^{\e}+vI_{[t,t+\e]}+\f^*,$
$X^*$, $X^\e$ are the state process on $[0,T]$ associated with $u^*$, $u^\e$, respectively.
\end{definition}

\medskip

We emphasize that both open-loop equilibrium strategy and closed-loop equilibrium strategy are independent of initial state $x_0$. However, the perturbed control $u^{v,\e}$ in Definition \ref{Definition-1} is actually different from $u^\e$ in Definition \ref{Definition-3}. In this paper, we will demonstrate further connections between these two kinds of strategies.

\medskip

In the following, let $K$ be a generic constant which varies in different context and
\bel{Simplified-notations-1}\ba{ll}
\ns\ds \sR:= R+\wt R,\ \ \sQ:=Q+\wt Q, \ \ \sG:=G+\wt G,\ \ \sS=S+\wt S.
\ea\ee

\section{Characterizations of equilibrium controls/strategies}

In this part, we state the main results of this article. We start with the case of open-loop equilibrium controls. To this end, given $u\in L^2_{\dbF}(0,T;\dbR^m)$, we introduce
\bel{Systems-for-open-only}\left\{
\ba{ll}
\ns\ds dP_1=-\Big[P_1 A+A^{\top}P_1+C^{\top}P_1C-Q\Big]ds,\ \ \\
\ns\ds dP_2=-\Big\{P_2A+A^{\top}P_2-\wt Q \Big\}ds,\\
\ns\ds d P_{3}=-\Big[A^{\top}P_3+P_2b+ (P_2B  -\wt S^{\top})u\Big]ds+ L_{3}dW(s),\\
d P_4=-\Big\{A^{\top} P_4+C^{\top} L_4+C^{\top} P_1\si+P_1 b+(C^{\top} P_1 D   \\
\ns\ds\qq \qq +P_1 B   -S^{\top})u\Big\}ds+L_4dW(s),\\
\ns\ds P_1(T)=-G,\ P_2(T)=-\wt G,\ P_3(T)=0,\ P_4(T)=-g.
\ea\right.\ee
Here $P_1, P_2$ do not rely on $u$ while $P_3, P_4$ do. It is easy to see the solvability, as well as the following regularities, of systems of equations (\ref{Systems-for-open-only}),
$$
\ba{ll}
\ns\ds P_1,\ P_2\in C([0,T];\dbR^{n\times n}),\ (P_{3 },\L_{3}),(P_{4},\L_{4})\in  L^2_{\dbF}(\Omega;C([0,T];\dbR^{n}))\times L^2_{\dbF}(0,T;\dbR^{n}).
\ea
$$
For $X$ in (\ref{state-equation-original}), we define
\bel{Arbitrary-s-Y-Z-general-2}\left\{
\ba{ll}
\ns\ds M(s,t):= P_1(s)X(s)+   P_2(s)\dbE_tX(s) + \dbE_t P_3(s)+ P_4(s), \ \ s\in[t,T],\\
\ns\ds N(s):= P_1(s)\big(C(s) X(s) +D(s)u(s) +\si(s)\big)+ L_4(s),\ \ s\in[0,T].
\ea\right.\ee
\bt\label{Equivalence-open-1}
Suppose (H1) holds, $P_1$ satisfies (\ref{Systems-for-open-only}). Then $\bar u $ is an open-loop equilibrium control associated with initial state $\bar X(0)=x_0\in\dbR^n$ if and only if
\bel{necessity-open-loop-1}\ba{ll}
\ns\ds \sR(s)-D(s)^{\top}P_1(s)D(s)\geq0,\qq s\in[0,T], \ \ a.e. \ \
\ea\ee
and given $(\bar M,\bar N)$ in (\ref{Arbitrary-s-Y-Z-general-2}) associated with $\bar u $,
\bel{necessary-open-loop-2}
\ba{ll}
\ns\ds  \sR(s) \bar u(s)+\sS(s)\bar X(s)- B(s)^{\top} \bar M(s,s)
- D(s)^{\top} \bar N(s)=0, \ \ s\in[0,T]. \ \ a.e.
\ea
\ee
\et

\medskip

Above (\ref{necessity-open-loop-1}), (\ref{necessary-open-loop-2}) are named as $first$-$order$, $second$-$order$ $equilibrium$ $conditions$, which are comparable with classical  $first$-$order$, $second$-$order$ $necessary$ $optimality$ $conditions$ (e.g., \cite{Chen-Li-Zhou-1998}, \cite{Yong-Zhou 1999}) in optimal control theories.

\br\label{Open-loop-second-order}
As to $P_1$ in (\ref{necessity-open-loop-1}), it is indeed the unique solution of classical second-order adjoint equation in optimal control theories. That is to say, (\ref{necessity-open-loop-1}) can reduce into the traditional second-order necessary optimality condition if $\wt R=0.$ To our best, this point was not discussed seriously in \cite{Hu-Jin-Zhou-2012}, \cite{Hu-Jin-Zhou-2017}, and other related papers on open-loop equilibrium controls.
\er
\br\label{Open-loop-first-order}
For $X$ in (\ref{state-equation-original}), we see that $(M,N)$ satisfies
\bel{Equations-for-M-N-1}\left\{\!\!\!\!\ba{ll}
\ns\ds dM=-\Big[A^{\top}M  + C^{\top} N-QX-S^{\top} u-\wt Q\dbE_t X-\wt S^{\top}\dbE_t u\Big]dr+NdW(r), \\
\ns\ds M(T,t)=-G X (T)-\wt G  \dbE_tX (T)-g.
\ea\right.\ee
As a result, if $\wt R=\wt Q=\wt S=\wt G=0$ and $u$ is optimal, (\ref{Equations-for-M-N-1}) becomes the first-order adjoint equation. In other words, (\ref{necessary-open-loop-2}) degenerates into an equivalent form of first-order necessary condition.
\er

\br
If $\wt R=\wt S=S=0,$ $R, Q, G$ are definite matrices, then (\ref{necessity-open-loop-1}) is obvious to see. In this scenario, a characterization of open-loop equilibrium control, which is different yet equivalent with (\ref{necessary-open-loop-2}), was given in Theorem 3.5 of \cite{Hu-Jin-Zhou-2017}. However, there were no systems of equations (\ref{Systems-for-open-only}) involved in their conclusion.
\er
Next we characterize the closed-loop representation of open-loop equivalent control in the sense of Definition \ref{Definition-2}. For $(\Th_2,\f)\in L^2(0,T;\dbR^{m\times n})\times L^2_{\dbF}(0,T;\dbR^m)$ in above (\ref{Definition-perturbed-introduction}), we introduce system of equations
\bel{System-closed-loop-representation}\left\{\!\!\!\!
\ba{ll}
\ns\ds d\cP_1\!=\!-\Big[\!\cP_1A\!+\!A^{\top}\cP_1\!+\!C^{\top}\cP_1C \!+\!(\cP_1B\!+\!C^{\top}\cP_1D\!-\!S^{\top})\Th_2\!-\!Q\Big]ds,\ \ \\
\ns\ds d\cP_2\!=\!-\Big\{\!\cP_2A \!+\!A^{\top}\cP_2\!-\! \wt Q\! +\!(\cP_2B\!-\!\wt S^{\top})\Th_2 ]\Big\}ds,\\
\ns\ds d\cP_{3}\!=\!-\Big[\!A^{\top}\cP_3\!+\! (\cP_2B\!-\!\wt S^{\top})\f\!+\!\cP_2b\Big]ds\!+\!\cL_{3}dW(s),\\
d\cP_4\!=\!-\Big\{\!A^{\top}\cP_4\!+\!C^{\top}\cL_4\!+\!C^{\top} \cP_1\si\!+\!
(C^{\top} \cP_1  D
\!+\!\cP_1 B\!-\!S^{\top})\f\\
\ns\ds\qq\qq +\!\cP_1b\!\Big\}ds\!+\!\cL_4dW(s),\\
\ns\ds \cP_1(T)=-G,\ \cP_2(T)=-\wt G,\ \cP_3(T)=0,\ \cP_4(T)=-g,
\ea\right.\ee
and following-up processes $(\cM,\cN)$ as follows,
\bel{Definition-c-M-N}\left\{\ba{ll}
\ns\ds \cM:=\cP_1X+\cP_2\dbE_t X+\dbE_t\cP_3+\cP_4, \\
\ns\ds \cN:=\cP_1(C+D\Th_2)X
+\cP_1(D\f+\si)+\cL_4.
\ea\right.\ee
\br
Given $(\Th_2,\f)$, if $u:=\Th_2 X+\f$ where $X$ is the associated state satisfying (\ref{state-equation-original}) on $[0,T]$, we see that $(\cM,\cN)$ solves (\ref{Equations-for-M-N-1}) as well. By the uniqueness of BSDEs, $(\cM,\cN)\equiv (M,N)$. Consequently, we obtain two different representations, i.e., (\ref{Arbitrary-s-Y-Z-general-2}), (\ref{Definition-c-M-N}), for the solutions of (\ref{Equations-for-M-N-1}).
\er
\begin{theorem}\label{Th-closed-loop-representation}
Suppose (H1) holds, $P_1$ satisfies (\ref{Systems-for-open-only}). Then for any $X^*(0)=x_0\in\dbR^n$, there exists equilibrium control $u^*$ in the sense of Definition \ref{Definition-2} if and only if (\ref{necessity-open-loop-1}) is true and there exist $ \cP_1^*,\ \cP_2^*,$ $(\cP_3^*,\cL^*_3)$, $(\cP_4^*,\cL^*_4)$ satisfying BSDEs (\ref{System-closed-loop-representation}) with $(\Th_2,\f)\equiv (\Th^*,\f^*)$ and
\bel{One-equality-for-Theta}\left\{\!\!\! \ba{ll}
\ns\ds \big[\sR- D^{\top}\cP_1^* D\big]\Th^*=B^{\top}\big[\cP_1^*+\cP_2^*\big]+D^{\top} \cP_1^* C-\sS,\\
\ns\ds
\big[\sR- D^{\top}\cP_1^* D\big]\f^*= D^{\top}\big[\cP_1^*\si+\cL_4^*\big]+ B^{\top}
\big[\cP_3^*+\cP_4^*\big].
\ea\right.\ee
\end{theorem}
\br
From (\ref{One-equality-for-Theta}), there exists $\th'\in L^2(0,T;\dbR^{n\times n})$, $\f'\in L^2_{\dbF}(0,T;\dbR^m) $ s.t.
\bel{representations-for-Th-f}\left\{\ba{ll}
\ns\ds \Th^*=\big[\sR-D^{\top}\cP_1^* D\big]^{\dagger}\big[B^{\top}(\cP_1^*+\cP_2^* )+D^{\top}\cP_1^* C-\sS\big]\\
\ns\ds\qq +\Big\{I-\big[\sR-D^{\top}\cP_1^*D\big]^{\dagger}\big[\sR-D^{\top}
\cP_1^*D\big]\Big\}\th',\\
\ns\ds \f^*=  \big[\sR-D^{\top}\cP_1^*D\big]^{\dagger}\big[B^{\top}[\cP_4^*+\cP_3^*]
+D^{\top}[\cP_1^*\si +\cL_4^*]\big]\\
\ns\ds \qq +\Big\{I-\big[\sR-D^{\top}\cP_1^*D\big]^{\dagger}\big[\sR-D^{\top}\cP_1^*D\big]\Big\}\f'.
\ea\right.\ee
Moreover,
\bel{Some-furthermore-condition-on}\left\{\ba{ll}
\ns\ds \cR\Big(B^{\top}(\cP_1^*+\cP_2^*)+D^{\top}\cP_1^* C-\sS\Big)\subset \cR\Big(\sR-D^{\top}\cP_1^* D\Big),\ \ a.e. \\
\ns\ds \Big[B^{\top}[\cP_4^*+\cP_3^*]
+D^{\top}[\cP_1^*\si +\cL_4^*]\Big]\in \cR\Big(\sR-D^{\top}\cP_1^* D\Big), \ a.e. \ a.s. \\
\ns\ds \big[\sR-D^{\top}\cP_1^* D\big]^{\dagger}\big[B^{\top}(\cP_1^*+\cP_2^*)+D^{\top}\cP_1^* C-\sS\big]\in L^2(0,T;\dbR^{m\times n}),\\
\ns\ds \Big[\big[\sR-D^{\top}\cP_1^*D\big]^{\dagger}\big[B^{\top}[\cP_4^*+\cP_3^*]
+D^{\top}[\cP_1^*\si +\cL_4^*]\big] \in L^2_{\dbF}(0,T;\dbR^m).
\ea\right.\ee
In above, $\cR(A)$, $A^{\dagger}$ is the range, pseudo-inverse of matrix $A$, respectively. Therefore, we obtain one representation of open-loop equilibrium strategy pair $(\Th^*,\f^*)$, as well as some intrinsic relations among coefficients in (\ref{Some-furthermore-condition-on}). Compared with open-loop equilibrium controls in Theorem \ref{Equivalence-open-1}, such closed-loop representations are advantageous in some sense and provide us more useful information.
\er

At last, we give the characterizations of closed-loop equilibrium strategies.
For $(\Th_1,\f)\in L^2(0,T;\dbR^{m\times n})\times L^2_{\dbF}(0,T;\dbR^m)$ in above (\ref{Definition-perturbed-introduction}), we introduce
\bel{System-of-closed-loop-strategy}\left\{
\ba{ll}
\ns\ds d\sP_1=-\Big[\sP_1 (A+B\Th_1)+(A+B\Th_1)^{\top}\sP_1+(C+D\Th_1)^{\top}\sP_1(C+D\Th_1)\\
\ns\ds\qq\qq -\big[Q+\Th_1^{\top}S+\Th_1^{\top}R\Th_1+S^{\top}\Th_1\big] \Big]ds,\ \ \\
\ns\ds d\sP_2=-\Big\{\sP_2(A+B\Th_1) +(A+B\Th_1)^{\top}\sP_2-\big[ \wt Q+\Th_1^{\top}\wt S+\Th_1^{\top}\wt R\Th_1+\wt S^{\top}\Th_1\big]\Big\}ds,\\
\ns\ds d\sP_{3}=-\Big[(A+B\Th_1)^{\top}\sP_3+ \sP_2b
+(\sP_2B -\wt S^{\top}-\Th_1^{\top}\wt R)\f \Big]ds+\sL_{3}dW(s),\\
d\sP_4=-\Big\{(A+B\Th_1)^{\top}\sP_4+(C+D\Th_1)^{\top}\sL_4+(C+D\Th_1)^{\top} \sP_1 (D\f+\si)\\
\ns\ds\qq \qq +\sP_1 (B\f+b)-(S^{\top}+\Th^{\top}_1R)\f\Big\}ds+\sL_4dW(s),\\
\ns\ds \sP_1(T)=-G,\ \sP_2(T)=-\wt G,\ \sP_3(T)=0,\ \sP_4(T)=-g,
\ea\right.\ee
and following-up $\sM,\ \sN$ as follows,
\bel{Definitions-of-s-M-s-N}\left\{\!\!\ba{ll}
\ns\ds \sM:=\sP_1 X+\sP_2\dbE_t X+\dbE_t\sP_3+\sP_4,\\
\ns\ds \sN:=\sP_1(C+D\Th_1)X+\sP_1(D\f+\si)+\sL_4.
\ea\right.\ee
\bt\label{Th-closed-loop-strategy}
A pair of $(\Th^*,\f^*)\in L^2(0,T;\dbR^{m\times n}) \times L^2_{\dbF}(0,T;\dbR^{m})$ is a closed-loop equilibrium strategy if and only if there exists $\sP_i^*$ satisfies (\ref{System-of-closed-loop-strategy}) with $(\Th_1,\f)\equiv (\Th^*,\f^*)$ such that
\bel{Conclusion-of-closed-strategy}\left\{\!\!\! \ba{ll}
\ns\ds \sR-D^{\top}\sP_1^* D\geq0,\\
\ns\ds (\sR-D^{\top}\sP_1^*D)\Th^*= B^{\top}(\sP_1^*+\sP_2^*)+D^{\top}\sP_1^*C-\sS,\\
\ns\ds (\sR-D^{\top}\sP_1^*D)\f^*=B^{\top}
(\sP_3^*+\sP_4^*)+D^{\top}\sP_1^* \si+D^{\top}\sL_4^*.
\ea\right.
\ee
\et

For the closed-loop equilibrium strategy $(\Th^*,\f^*)$, the first inequality in (\ref{Conclusion-of-closed-strategy}) is referred as the $second$-$order$ $equilibrium$ $condition$, while the other two conditions are named as $first$-$order$ $equilibrium$ $condition$.

\br
We make some comparisons among (\ref{Systems-for-open-only}), (\ref{System-closed-loop-representation}), (\ref{System-of-closed-loop-strategy}), from which we see the connections between open-loop equilibrium controls and their closed-loop representations, as well as that of closed-loop equilibrium controls and closed-loop representations.

$\diamond$ The later two systems reduce to the first one if $\Th_1=0$, or $\Th_2=0$, and $\f\equiv u$.

$\diamond$ The solutions of the first two equations in (\ref{Systems-for-open-only}), (\ref{System-of-closed-loop-strategy}) are symmetric, while the analogue of (\ref{System-closed-loop-representation}) are non-symmetric (see e.g., \cite{Yong-2017}).

$\diamond$ The first two equations in (\ref{Systems-for-open-only}) merely depends on given coefficients, while the counterparts in (\ref{System-closed-loop-representation}) and (\ref{System-of-closed-loop-strategy}) are determined by $\Th_1$, or $\Th_2$.

$\diamond$ The last two equations in (\ref{Systems-for-open-only}) rely on control process $u$, while the analogue equations in (\ref{System-closed-loop-representation}) and (\ref{System-of-closed-loop-strategy}) are determined by $\f$.
\er

\medskip

\br\label{Relation-between-open-loop-closed-loop}
To capture the new feature of time inconsistency, let $\wt G=\wt S=\wt Q=\wt R=0$, $b=\si=g=0$. Suppose there exists closed-loop representation of open-loop optimal control $u_1^*=\Th_1^* X_1^*$ and closed-loop optimal control $u_2^*:=\Th_2^* X_2^*$, where $\f_1^*=\f_2^*=0$. We claim that
$
\sP_1^*=\cP_1^*.$
If furthermore $R-D^{\top}\cP_1^*D>0$, a.e., $u^*_1=u_2^*$, and $\cP_1^*\equiv\sP_1^*$ satisfies the Riccati equations in classical stochastic linear quadratic problems.
Actually, in this setting,
$$\cP_2^*=\cP_3^*=\cL_3^*=\cP_4^*=\cL_4^*=0,\ \ \sP_2^*=\sP_3^*=\sL_3^*=\sL_4^*=\sP_4^*=0, \ $$
and the last two conditions in (\ref{One-equality-for-Theta}),  (\ref{Conclusion-of-closed-strategy}) become,
\bel{Combined-two-first-order-special}\ba{ll}
\ns\ds
\big[R- D^{\top}\cP_1^* D\big]\Th_1^*=B^{\top}\cP_1^*+D^{\top} \cP_1^* C-S,\\
\ns\ds
(R-D^{\top}\sP_1^*D)\Th_2^*= B^{\top}\sP_1^*+D^{\top}\sP_1^*C-S.
\ea
\ee
Substituting the second expression into the first equation of (\ref{System-of-closed-loop-strategy}) with $(\Th_1,\f)\equiv (\Th^*_2,0)$, we have
$$\ba{ll}
\ns\ds d\sP_1^*=-\Big[\sP_1^*A+A^{\top}\sP_1^*+C^{\top}\sP_1^*C +(\sP_1^*B+C^{\top}\sP_1^*D-S^{\top})\Th^*_2-Q\Big]ds.
\ea$$
For $u\in L^2_{\dbF}(0,T;\dbR^m)$, $\xi\in L^2_{\cF_t}(\Omega;\dbR^n)$, by It\^{o}'s formula to $X^{\top}\sP_1^* X$,
$$\ba{ll}
\ns\ds   X(T)^{\top}GX(T) + \int_t^T\big[X^{\top}QX+2u^{\top}SX
+u^{\top}Ru\big]dr\\
\ns\ds =-\xi^{\top}\sP^*_1(t)\xi+\int_t^T\cL_2dr+\int_t^T\cW_2 dW(r),
\ea
$$
where
$$\left\{\ba{ll}
\ns\ds \cL_2:=u^{\top}\h R_2 u+2u^{\top}\h S_2X-X^{\top}\h S_2^{\top}\Th_2^* X, \\
\ns\ds\cW_2:=-2X^{\top}\big[\sP_1^*CX+\sP_1^* Du\big],\\
\ns\ds \h S_2:=S-D^{\top}\sP_1^*C-B^{\top}\sP_1^*,\ \ \h R_2:=R-D^{\top}\sP_1^*D.
\ea\right.
$$
Thanks to the second equality of (\ref{Combined-two-first-order-special}), as well as the symmetry of $R,\ \sP_1^*$,
$$
-X^{\top}\h S^{\top}_2\Th_2^*X=X^{\top}\big[\Th_2^*\big]^{\top}\h R_2 \Th_2^* X,\ \ u^{\top}\h S_2 X=-u^{\top}\h R_2\Th^*_2X.
$$
As a result,
$$\cL_2=(u-\Th_2^*X)^{\top}\h R_2 (u-\Th_2^* X),$$
and for optimal control $u^*_2$, one has,
$$\ba{ll}
\ns\ds  V(t,\xi)=J(t,\xi,u^*_2)=-\frac  1 2
\xi^{\top}\sP_1^*(t)\xi.
\ea
$$
Similarly we can deduce that
$$\ba{ll}
\ns\ds V(t,\xi)=J(t,\xi,u^*_1)=-\frac  1 2
\xi^{\top}\cP_1^*(t)\xi.
\ea
$$
By the continuity of $\cP_1^*,\sP_1^*,$ and the arbitrariness of $\xi$,
$$\ba{ll}
\ns\ds \dbP\big\{\omega\in\Omega;\ \cP_1^*(t,\omega)=\sP_1^*(t,\omega),\ \ \forall t\in[0,T]\big\}=1.
\ea
$$
The equality of $u^*_1=u_2^*$ is easy to obtain.

In general, $\cP_1^*$ is different from $\sP_1^*$, not to mention the equality of $u^*_1=u_2^*$. For example, when $\wt G\neq0$, one can see that $\sP_1^*$ is symmetric while $\cP_1^*$ is not.

To sum up, the closed-loop optimal controls coincide with closed-loop representation of open-loop optimal controls under proper conditions. However, this relation breaks when time-inconsistency happens.
\er
\br
For the second-order equilibrium conditions in Theorem \ref{Equivalence-open-1}, Theorem \ref{Th-closed-loop-representation} and Theorem \ref{Th-closed-loop-strategy}, we have the following comments.

$\diamond$ As to open-loop equilibrium controls, no matter it has closed-loop representations or not, we use $\sR-D^{\top}\cP_1^* D\geq0$, where $\cP_1^*$ satisfies the second-order adjoint equation in classical stochastic linear quadratic optimal control problems. This condition was missing in \cite{Hu-Jin-Zhou-2012}, \cite{Hu-Jin-Zhou-2017}, \cite{Wei-Wang-2017}, \cite{Wang-Wu-2015}.

$\diamond$ As to closed-loop equilibrium controls, we introduce $\sR-D^{\top}\sP_1^* D\geq0$ where $\sP_1^*$ satisfies one backward ordinary differential equation that contains Riccati equation as special case. Notice that this condition has not been discussed in \cite{Bjork-Murgoci-2012}, \cite{Bjork-Murgoci-2014}, \cite{Yong-2012-2}, \cite{Yong-2017}.
\er

\br\label{Relation-between-open-loop-closed-loop-2}
At this moment, we revisit the open-loop equilibrium controls and closed-loop equilibrium controls when $\wt G=\wt S=\wt Q=\wt R=g=0$.

From Remark \ref{Open-loop-second-order}, \ref{Open-loop-first-order}, the open-loop equilibrium controls under this framework are fully characterized by first-order, second-order necessary optimality conditions. This gives us a quantitative and clear picture of this kind of equilibrium control. Notice that the characterization of optimal controls includes first-order necessary condition and the following convexity condition (see \cite{Chen-Yong-2001})
\bel{convexity-condition}\ba{ll}
\ns\ds \dbE_t\int_t^T u^{\top}\big[R u+S X^{0}-B^{\top}Y^{0}-D^{\top}Z^{0}\big]dr\geq 0,\ \ \ \forall u\in L^2_{\dbF}(t,T;\dbR^m),
\ea
\ee
where $X^0$ satisfies (\ref{state-equation-original}) with $\xi=0$, $(Y^0,Z^0)$ solves (\ref{Equations-for-M-N-1}) with $\wt G=\wt S=\wt Q=g=0$ and $X\equiv X^0$.
Consequently, the exact difference between equilibrium controls and optimal controls in the open-loop sense is attributed to that between (\ref{necessity-open-loop-1}) and (\ref{convexity-condition}).

For closed-loop equilibrium controls/strategies in Theorem \ref{Th-closed-loop-strategy}, their characterization (\ref{Conclusion-of-closed-strategy}) reduces to
\bel{Conclusion-of-closed-strategy-classical}\ba{ll}
\ns\ds R-D^{\top}\sP_1^* D\geq0,\ \ \ \ (R-D^{\top}\sP_1^*D)\Th^*= B^{\top}\sP_1^*+D^{\top}\sP_1^*C-S,\\
\ns\ds (R-D^{\top}\sP_1^*D)\f^*=B^{\top}
\sP_4^*+D^{\top}\sP_1^* \si+D^{\top}\sL_4^*.
\ea
\ee
According to \cite{Sun-Li-Yong-2016}, \cite{Sun-Yong-2014}, (\ref{Conclusion-of-closed-strategy-classical}) is equivalent to the optimality of strategy pair $(\Th^*,\f^*)$ or control variable $u^*:=\Th^*X^*+\f^*.$ In other words, our defined closed-loop equilibrium controls/strategies are natural extension of closed-loop optimal controls/strategies. This not only leads to one more essential distinction between open-loop, closed-loop equilibrium controls, but not illustrate the reasonability of introduced closed-loop equilibrium controls from the optimality viewpoint.
\er

\section{Proofs of the main results}

In this section, we prove Theorem \ref{Equivalence-open-1}--\ref{Th-closed-loop-strategy}.

For $(\Th_1,\Th_2,\f)\in L^2(0,T;\dbR^{m\times m})\times L^2(0,T;\dbR^{m\times m})\times  L^2_{\dbF} (0,T;\dbR^{m})$, we consider
\bel{State-equation-Th-1-Th-2}\left\{\2n\ba{ll}
\ns\ds
dX =\big[A  X +B (\Th_1+\Th_2) X +B\f +b \big]ds \\
\ns\ds\qq \q  +\big[C X
+D  (\Th_1+\Th_2) X+D \f
+\si  \big]dW(s), \ \ s\in[0,T],\\
\ns\ds  X(0)=x_0.\ea\right.\ee
In the following, let
\bel{Definition-of-u-any}\ba{ll}
\ns\ds u:= (\Th_1+\Th_2)X+\f, \ \ u^{\e}:=\Th_1 X^{\e}+\Th_2 X+\f+vI_{[t,t+\e]}.
\ea\ee
Fix $t\in[0,T)$, $v\in\dbR^m$ and small $\e>0$, let $ X^\e $ be the solution to the following perturbed system:
\bel{Perturbed-1.1-ep}\left\{\2n\ba{ll}
\ns\ds
dX^\e =\big[(A+B\Th_1) X^\e+B\Th_2 X+B vI_{[t,t+\e]}  + B\f
+b \big]ds\\
\ns\ds\qq\q  +\big[(C+D \Th _1) X^\e +D\Th_2 X+ D vI_{[t,t+\e]}
 +D\f+\si \big]dW(s),   \\
\ns\ds X^\e(0)=x_0,
\ea\right.\ee
with $s\in[0,T].$ Hence we see that $X_0^\e:=X^\e-X$ satisfies
\bel{Equation-of-s-X-1-epislon}\left\{\2n\ba{ll}
\ns\ds
dX_0^\e =\big[(A +B \Th_1) X^{\e}_0   + B vI_{[t,t+\e]}  \big]ds\\
\ns\ds \qq\qq +\big[(C +D \Th_1) X^{\e}_0 + D vI_{[t,t+\e]}\big]dW(s), \\
\ns\ds X^\e_0(0)=0.\ea\right.\ee
\br\label{estimates}
By Proposition 2.1 in \cite{Sun-Yong-2014},  we have the following estimate of $X_0^\e$
$$\ba{ll}
\ns\ds \dbE_t\sup_{r\in[t,t+\e]}|X^\e_0(r)|^2\leq K\e,\ \  a.s., \ \ t\in[0,T).
\ea
$$
\er
To begin with, we have the following difference of cost functional.

\medskip

\bl\label{Cost-functional-1}
 Suppose (H1) holds, $(\Th_1,\Th_2,\f)$ are given as above, $u$, $u^\e$ are defined in (\ref{Definition-of-u-any}). Then we have
\bel{difference-cost}\ba{ll}
\ns\ds J(t,x,u^\e(\cd))-J(t,x,u(\cd)) =
J_1(t,x)+J_2(t,x)+\dbE_t\int_t^{t+\e}\lan (\sS^{\top}+\Th_1^{\top}\sR)v,X_0^\e\ran ds,
\ea\ee
where $\sR,\ \sS$ are defined in (\ref{Simplified-notations-1}),
$$\left\{\!\!\!\ba{ll}
\ns\ds J_1(t):=\dbE_t \int_t^T \big[\lan F_1,X_0^\e\ran +\1n\lan F_2,vI_{[t,t+\e)}\ran\big]ds +\dbE_t\lan GX(T)+\wt
G \dbE_tX(T)+g,X_0^\e(T)\ran,\\
\ns\ds J_2(t):= {1\over2}\dbE_t\int_t^T\2n \lan F_1^\e, X_0^\e \ran  ds+ \frac 1 2 \dbE_t \lan GX_0^\e(T)+ \wt G \dbE_tX_0^\e(T),X_0^\e(T)\ran,
\ea\right.
$$
and
$$\left\{\ba{ll}
\ns\ds F_1 \equiv \big[Q+\Th_1^{\top}S+\Th_1^{\top}R(\Th_1+\Th_2)+S^{\top}(\Th_1+\Th_2)\big] X +(S^{\top}+\Th^{\top}_1R)\f\\
\ns\ds\qq +\big[\wt Q+\Th_1^{\top}\wt S+\Th_1^{\top}\wt R(\Th_1+\Th_2)+\wt S^{\top}(\Th_1+\Th_2)\big]\dbE_t X  +(\wt S^{\top}+\Th^{\top}_1\wt R)\dbE_t \f,\\
\ns\ds F_2 \equiv \frac 1 2 \sR v+\big[ S +R(\Th_1+\Th_2)\big] X +R\f
+\big[\wt S +\wt R(\Th_1+\Th_2)\big]\dbE_t X+\wt R\dbE_t\f,\\
\ns\ds F_1^\e\equiv \big[Q + S^{\top}\Th_1 +\Th_1^{\top}S+\Th_1^{\top}R\Th_1\big]X_0^\e
+\big[\wt Q +\wt S^{\top}\Th_1+\Th_1^{\top}\wt S+\Th_1^{\top}\wt R\Th_1\big]\dbE_t X_0^\e.
\ea\right.
$$
\el

\begin{proof} By above definitions of $X$, $X^{\e}$ and $X_0^{\e}$, we deal with the terms in the cost functional one by one. First let us treat the term associated with $Q$,
$$\ba{ll}
\ns\ds \lan QX^\e,X^{\e}\ran-\lan QX,X\ran=2\lan QX,X_0^\e\ran+\lan QX_0^\e,X_0^\e\ran.
\ea$$
Then we look at the one with $S$. From the definitions of $u$ and $u^\e$, we have
$$\ba{ll}
\ns\ds   \lan S X^\e ,u^{\e} \ran- \lan S X ,u \ran\\
\ns\ds= \lan S^{\top}\Th_1 X_0^\e,X_0^\e\ran+ \lan X_0^\e, S^{\top}
\big[(\Th_1+\Th_2) X+vI_{[t,t+\e]} +\f \big]\ran \\
\ns\ds\q + \lan X_0^\e,\Th_1^{\top}SX\ran+ \lan SX, vI_{[t,t+\e]}\ran.
\ea$$
We also have
$$ \ba{ll}
\ns\ds  \lan R u^\e,u^{\e}\ran-\lan Ru,u\ran\\
\ns\ds=\lan \Th_1^{\top}R\Th_1 X_0^\e, X_0^\e\ran+ 2\lan RvI_{[t,t+\e]},\Th _1X_0^\e\ran+\lan Rv,vI_{[t,t+\e]}\ran \\
\ns\ds\q+2 \lan R\Th_1 X_0^\e,(\Th_1+\Th_2)X+\f\ran
+2\lan R vI_{[t,t+\e]},(\Th_1+\Th_2)X+\f\ran.
\ea $$
Similarly one can obtain the terms involving $\wt Q$, $\wt S$, $\wt R$ as,
$$\left\{\ba{ll}
\ns\ds \lan \wt Q \dbE_tX^\e,\dbE_tX^{\e}\ran-\lan \wt Q\dbE_tX,\dbE_tX \ran =2\lan \wt Q \dbE_tX ,\dbE_t X_0^\e\ran+\lan \wt Q\dbE_t X_0^\e,\dbE_t X_0^\e \ran,\\
\ns\ds   \lan \wt S \dbE_tX^\e ,\dbE_tu^{\e} \ran- \lan \wt S \dbE_tX ,\dbE_tu \ran\\
%
%
\ns\ds= \lan \wt S^{\top}\Th_1 \dbE_tX_0^\e,\dbE_tX_0^\e\ran+ \lan \dbE_tX_0^\e, \wt S^{\top}
\big[(\Th_1+\Th_2) \dbE_tX+vI_{[t,t+\e]} +\dbE_t\f \big]\ran\\
\ns\ds\q + \lan \dbE_tX_0^\e,\Th_1^{\top}\wt S\dbE_tX\ran+ \lan \wt S\dbE_tX, vI_{[t,t+\e]}\ran,\\
\ns\ds  \lan \wt R \dbE_tu^\e,\dbE_tu^{\e}\ran-\lan \wt R\dbE_tu,\dbE_tu\ran\\
\ns\ds = \lan \Th_1^{\top}\wt R\Th_1 \dbE_tX_0^\e, \dbE_tX_0^\e\ran+ 2\lan \wt RvI_{[t,t+\e]},\Th _1\dbE_tX_0^\e\ran+\lan \wt Rv,vI_{[t,t+\e]}\ran \\
\ns\ds\q+2 \lan \wt R\Th_1 \dbE_tX_0^\e,(\Th_1+\Th_2)\dbE_tX+\dbE_t\f\ran
+2\lan \wt R vI_{[t,t+\e]},(\Th_1+\Th_2)\dbE_tX+\dbE_t\f\ran.
\ea\right.$$
At last we have the follows results on the terms associated with $ G $ and $\wt G $,
$$\left\{\ba{ll}
\ns\ds \lan G X^\e(T),X^\e(T)\ran-\lan G X(T),X(T)\ran\\
\ns\ds=2\lan G X(T),X_0^\e(T)\ran+\lan G X_0^\e(T),X_0^\e(T)\ran,\\
\ns\ds \lan \wt  G \dbE_tX^\e(T),\dbE_tX^\e(T)\ran-\lan \wt G \dbE_tX(T),\dbE_tX(T)\ran\\
\ns\ds=2\lan \wt  G \dbE_tX(T),\dbE_t X_0^\e(T)\ran+\lan \wt G \dbE_tX_0^\e(T),\dbE_tX_0^\e(T)\ran.
\ea\right.$$
To sum up, we deduce above (\ref{difference-cost}).
\end{proof}

Next we spread out further study on $J_1(t)$ and $J_2(t)$ by making some equivalent transformations. In fact, from the definitions of equilibrium controls it is unavoidable to take certain convergence arguments. Fortunately, in above we derive
the important and useful structure of $\dbE_t\int_t^{t+\e}\lan F_2(r),v\ran dr$.
Consequently, we will derive similar expressions for other terms in $J_1(t)$, $J_2(t)$. This is the starting point for our later investigations.

\subsection{A new decoupling result}

Inspired by the decoupling tricks in the literature (e.g., \cite{Hu-Jin-Zhou-2012}, \cite{Yong-2013-SICON}, etc), we present one conclusion which serves our purpose of this paper. It is interesting in its own right and may be potentially useful for (among others) various problems.

Given $t\in[0,T]$, we consider
\bel{System-of-equations-1}\left\{\ba{ll}
\ns\ds dX=\big[A_1 X+A_2 \big]dr +\big[B_1 X + B_2 \big]dW(r),\ \ r\in[t,T],\\
\ns\ds dY=-\Big[C_1 Y+C_2 Z +C_3 X+C_4\dbE_t X+C_5 +\dbE_t C_6\Big]dr+ZdW(r),\\
\ns\ds X(0)=x,\ \   Y(T,t)=D_1 X(T)+D_2 \dbE_tX(T)+D_3.
\ea\right.
\ee

\mds

(H1) For $H:=\dbR^{m},\ \dbR^{n},\ \dbR^{n\times n}$, etc, suppose
$A_1,\ B_1,\ C_i\in L^2(0,T;H),$ $A_2,\ C_5\in L^2(\Omega;L^1(0,T;H)),$ $B_2\in L^2_{\dbF}(0,T;H),$ $D_1,\ D_2,\ D_3,\ x\in H.$

\bigskip

For $t\in[0,T]$ and $s\in[t,T]$, suppose that
\bel{Representation-for-Y-*}\ba{ll}
\ns\ds  Y(s,t)=P_1(s)X(s)+ P_2(s)\dbE_tX(s)+\dbE_t P_{3}(s)+ P_{4}(s),
\ea
\ee
where $P_1,\ P_2$ are deterministic, $P_3, \ P_4$ are stochastic processes satisfying
$$\left\{
\ba{ll}
\ns\ds  d P_i(s)=\Pi_{i}(s)ds,\ \ i=1,2, \ \ P_1(T)=D_1,\ \  P_2(T)=D_2, \\
\ns\ds d P_j(s)=\Pi_j(s)ds+\cL_j(s)dW(s),\ \ j=3,4,\ \  P_3(T)=0,\ \ P_{4}(T)=D_3.
\ea\right.
$$
Here $\Pi_i$ are to be determined. It is easy to see
$$\ba{ll}
\ns\ds d\dbE_t X=\big[A_1 \dbE_t X  +\dbE_tA_2 )\big]dr.
\ea
$$
Using It\^{o}'s formula, we derive that
$$\left\{
\ba{ll}
\ns\ds d\big[  P_1 X\big]=\Big[\Pi_{1} X+ P_1(A_1 X +A_2)\Big]ds
+ P_1\big(B_1 X  +B_2\big)dW(s),\\
\ns\ds d\big[P_2\dbE_t X\big]=\Big\{\Pi_2 \dbE_t X+P_2
\big[A_1\dbE_t X +\dbE_t A_2 \big]\Big\}ds.
\ea\right.
$$
As a result, we have
$$
\ba{ll}
\ns\ds d Y =\Big\{\big[\Pi_{1}+ P_1A_1\big]X+( \Pi_2+P_2A_1 )\dbE_tX
\\
\ns\ds\qq +\dbE_t\big[ \Pi_3+P_2A_2\big] +\Pi_4+P_1A_2\Big\}ds+\Big[ P_1B_1 X+P_1B_2+L_{4}\Big]dW(s).
\ea
$$
Consequently, it is necessary to see
\bel{Representation-for-Z-*}
\ba{ll}
\ns\ds Z=P_1B_1 X +P_1B_2+L_{4}.
\ea
\ee
In this case, from (\ref{Representation-for-Y-*}), (\ref{Representation-for-Z-*}), we see that
$$\left\{
\ba{ll}
\ns\ds \dbE_t Y=(P_1+P_2)\dbE_tX +\dbE_t \big[P_{3}+P_{4}\big],\\
\ns\ds \dbE_t Z=P_1B_1\dbE_tX+\dbE_t\big[P_1B_2 +L_{4}\big].
\ea\right.
$$
On the other hand,
$$
\ba{ll}
\ns\ds -\Big[C_1 Y+C_2 Z+C_3 X+C_4\dbE_t X+C_5+\dbE_t C_6 \Big]\\
\ns\ds=  -C_1\Big\{P_1X+P_2\dbE_t X +\dbE_tP_{3}+P_{4}\Big\} -C_2\Big[P_1B_1 X +P_1B_2+L_{4} \Big]\\
\ns\ds\qq -C_3 X-C_4\dbE_t X-C_5-\dbE_t C_6.
\ea
$$
At this moment, we can choose $\Pi_i(\cd)$ in the following ways,
$$\left\{
\ba{ll}
\ns\ds 0=\Pi_1+P_1 A_1+C_1P_1+C_2P_1B_1+C_3,\ \    \\
\ns\ds 0=\Pi_2 +P_2A_1 +C_1P_2+C_4,
\\
\ns\ds 0=\Pi_{4} +P_1  A_2  +C_1P_4+C_2\big[P_1 B_2 +L_4\big]+C_5,\\
\ns\ds 0=\Pi_3  +P_2  A_2  +C_1P_3+C_6.
\ea\right.
$$
Next we make above arguments rigorous.
Given the notations in (\ref{Simplified-notations-1}), for $s\in[0,T]$, we consider the following systems of equations
\bel{General-system-of-equations-needed}\left\{
\ba{ll}
\ns\ds dP_1=-\Big[P_1 A_1+C_1P_1+C_3P_1B_1+C_3\Big]ds,\ \ \\
\ns\ds dP_2=-\Big\{P_2A_1 +C_1P_2+C_4\Big\}ds,\\
\ns\ds dP_{3}=-\Big[C_1P_3+ P_2  A_2 +C_6 \Big]ds+L_{3}dW(s),\\
dP_4=-\Big\{C_1P_4+C_2L_4+C_2 P_1 B_2+P_1  A_2 +C_5\Big\}ds+L_4dW(s),\\
\ns\ds P_1(T)=D_1,\ P_2(T)=D_2,\ P_3(T)=0,\ P_4(T)=D_3.
\ea\right.\ee
From Proposition 2.1 in \cite{Sun-Yong-2014}, under (H1) we see the following regularities,
$$
\ba{ll}
\ns\ds P_1,\ P_2\in C([0,T];\dbR^{n\times n}),\ (P_{3},L_{3}),(P_{4},L_{4})\in  L^2_{\dbF}(\Omega;C([0,T];\dbR^{n}))\times L^2_{\dbF}(0,T;\dbR^{n}).
\ea
$$
At this moment, for $s\in[0,T]$, and $t\in[0,s]$, we define a pair of processes
\bel{Definition-s-Y-s-Z-1}
\ba{ll}
\ns\ds M:=P_1 X +P_2\dbE_tX+\dbE_tP_3 +P_4,  \ \  N:=P_1B_1 X  +P_1B_2+L_{4}.
%
\ea \ee
By the results of $P_i$, we can conclude that
$$(M_d,N)\in L^2_{\dbF}(\Omega;C([0,T];\dbR^n))\times L^2_{\dbF}(0,T;\dbR^n)$$
where $M_{d}(s)\equiv M(s,s)$ with $s\in[0,T]$.
We present the following result.

\bigskip

\bl\label{representation-for-Y-Z-*}
\rm Given $(\Th,\f)\in L^2(0,T;\dbR^{m\times n})\times L^2_{\dbF}(0,T;\dbR^m)$, suppose $(X,Y,Z)$ is the unique solution of (\ref{System-of-equations-1}) and $(M,N)$ are defined in (\ref{Definition-s-Y-s-Z-1}). Then for any $t\in[0,T]$,
$$\ba{ll}
\ns\ds   \dbP\Big\{\omega\in\Omega;\ Y(s,t)=M(s,t),\ \ \forall s\in[t,T] \Big\}=1,\\
\ns\ds \dbP\Big\{\omega\in\Omega;\ Z (s,t) =N(s)\Big\}=1,\ \ s\in[t,T]. \ \ a.e.
\ea
$$
\el
\begin{proof}
 Given (\ref{Definition-s-Y-s-Z-1}), it is easy to see that
$$
\ba{ll}
\ns\ds \dbE_tM =(P_1+P_2)\dbE_tX+  \dbE_t[P_3+P_4],\ \  \dbE_tN =P_1B_1\dbE_tX  +P_1\dbE_tB_2 +\dbE_t L_{4}.
\ea $$
Using It\^{o}'s formula, we know that
$$\left\{
\ba{ll}
\ns\ds d\big[P_1 X\big]=\Big[-(C_1P_1+C_2P_1B_1+C_3) X+P_1 A_2 \Big]ds +
P_1\big(B_1 X+B_2\big)dW(s),\\
\ns\ds d\big[P_2\dbE_t X\big]=\Big\{-\Big[C_1P_2 +C_4 \Big] \dbE_t X
+P_2\dbE_tA_2\Big\}ds.
\ea\right.
$$
Consequently, after some calculations one has
$$
\ba{ll}
\ns\ds dM=-\Big[C_1M +C_2N+C_3 X+C_4\dbE_t X+C_5 +\dbE_t C_6\Big]dr+NdW(r).
\ea$$
Considering $P_i(T)$ in (\ref{General-system-of-equations-needed}), we see that for any $t\in[0,T]$, $(M,N)\in L^2_{\dbF}(\Omega;C([t,T];\dbR^n))\times L^2_{\dbF}(0,T;\dbR^n)$ satisfies the backward equation in (\ref{System-of-equations-1}). The conclusion is followed by the uniqueness of BSDEs.
\end{proof}

\subsection{A new expression of $J_1$}

In this part, we deal with $J_1(t)$ in Lemma \ref{Cost-functional-1}. For convenience, we rewrite the equation of $X_0^\e:=X^\e-X$ as
\bel{Rewrite-of-X-0}\left\{\2n\ba{ll}
\ns\ds
dX_0^\e =\big[A_\th X^{\e}_0   + B vI_{[t,t+\e]}  \big]ds +\big[C_\th X^{\e}_0 + D vI_{[t,t+\e]}\big]dW(s), \\
\ns\ds X^\e_0(0)=0,\ea\right.\ee
where $s\in[0,T],$ and
$$
A_\th:= A+B\Th_1,\ \ \  C_\th:=C+D\Th_1.
$$
We introduce
\bel{Equations-for-Y-Z}\left\{\ba{ll}
\ns\ds dY=-\Big[A_\th^{\top}Y  + C_\th^{\top} Z-F_1\Big]dr+ZdW(r),\ \ r\in[t,T],\\
\ns\ds Y(T,t)=-G X (T)-\wt G  \dbE_tX (T)-g,
\ea\right.\ee
where $X$ satisfies (\ref{State-equation-Th-1-Th-2}), $F_1$ is in Lemma \ref{Cost-functional-1}. From Proposition 2.1 in \cite{Sun-Yong-2014}, (\ref{Equations-for-Y-Z}) is solvable with
$$
\ba{ll}
\ns\ds  (Y,Z) \in  L^2_{\dbF}(\Omega;C([t,T];\dbR^{n}))\times L^2_{\dbF}(t,T;\dbR^{n}),\ \ t\in[0,T).
\ea
$$
By It\^{o}'s formula on $[t,T]$, we have
$$\ba{ll}
\ns\ds d\lan Y,X_0^\e \ran=-\lan  A_\th^{\top} Y + C_\th^{\top} Z -F_1, X_0^\e \ran dr+\lan Z, X_0^\e \ran  dW(r)\\
\ns\ds\qq\qq\q+\lan Y, A_\th X_0^\e+B vI_{[t,t+\e]} \ran dr+\lan Y, C_\th X_0^\e + D vI_{[t,t+\e]}\ran dW(r)\\
\ns\ds\qq \qq\q+\lan Z, C_\th X_0^\e + D vI_{[t,t+\e]}  \ran dr.
\ea
$$
From (\ref{Equations-for-Y-Z}) we then arrive at
\bel{First-one-to-J-1}\ba{ll}
\ns\ds \dbE_t\lan -GX (T)-\wt G \dbE_tX (T)-g,X_0^\e(T)\ran-\dbE_t\int_t^T \lan F_1,X_0^\e\ran dr\\
\ns\ds =\dbE_t\int_t^{t+\e}\lan B^{\top} Y+ D^{\top} Z , v\ran dr.
\ea
\ee
Inspired by Lemma \ref{representation-for-Y-Z-*}, we introduce
\bel{Riccati-equation-system}\left\{\!\!
\ba{ll}
\ns\ds d\cP_1=-\Big[\cP_1 (A+B\Th_1+B\Th_2)+(C+D\Th_1)^{\top}\cP_1(C+D\Th_1+D\Th_2)\\
\ns\ds\qq\qq +(A+B\Th_1)^{\top}\cP_1-\big[Q+\Th_1^{\top}S+\Th_1^{\top}R(\Th_1+\Th_2)+S^{\top}(\Th_1+\Th_2)\big] \Big]ds,\ \ \\
\ns\ds d\cP_2=-\Big\{\cP_2(A+B\Th_1+B\Th_2) +(A+B\Th_1)^{\top}\cP_2-\big[ \wt Q+\Th_1^{\top}\wt S\\
\ns\ds\qq \qq +\Th_1^{\top}\wt R(\Th_1+\Th_2)+\wt S^{\top}(\Th_1+\Th_2)\big]\Big\}ds,\\
\ns\ds d\cP_{3}=-\Big[(A+B\Th_1)^{\top}\cP_3+ \cP_2(B\f+b)-(\wt S^{\top}+\Th_1^{\top}\wt R)\f \Big]ds+\cL_{3}dW(s),\\
d\cP_4=-\Big\{(A+B\Th_1)^{\top}\cP_4+(C+D\Th_1)^{\top}\cL_4+(C+D\Th_1)^{\top} \cP_1 (D\f+\si)\\
\ns\ds\qq \qq +\cP_1 (B\f+b)-(S^{\top}+\Th^{\top}_1R)\f\Big\}ds+\cL_4dW(s),\\
\ns\ds \cP_1(T)=-G,\ \cP_2(T)=-\wt G,\ \cP_3(T)=0,\ \cP_4(T)=-g.
\ea\right.\ee
Moreover, the following equalities hold on $[t,T],$
$$
\ba{ll}
\ns\ds Y=\cP_1 X +\cP_2 \dbE_tX  + \dbE_t\cP_3 +\cP_4,  \ \ Z= \cP_1(C+D\Th_1+D\Th_2) X  +\cP_1(D\f+\si)+\cL_{4}.
%
\ea $$
Consequently,
$$\ba{ll}
\ns\ds B^{\top}Y+D^{\top}Z=\big[B^{\top}\cP_1+D^{\top}\cP_1(C+D\Th_1+D\Th_2)\big]X
+B^{\top}\cP_2\dbE_t X\\
\ns\ds\qq \qq\qq\q  +B^{\top}\dbE_t \cP_3+B^{\top}\cP_4+D^{\top} \cP_1(D\f+\si)+D^{\top} \cL_4.
\ea
$$
This shows that
$$\ba{ll}
\ns\ds \dbE_t\int_t^{t+\e}\lan B^{\top} Y+ D^{\top} Z , v\ran dr\\
\ns\ds =\dbE_t\int_{t}^{t+\e}\lan \big[B^{\top}(\cP_1+\cP_2)+D^{\top}\cP_1(C+D\Th_1+D\Th_2)\big]X\\
\ns\ds\qq\qq\q +B^{\top}
(\cP_3+\cP_4)+D^{\top}\cP_1(D\f+\si)+D^{\top}\cL_4,v\ran dr.
\ea
$$
By the definition of $J_1(t)$ and above (\ref{First-one-to-J-1}), we see that
\bel{Expression-of-J-1}
\ba{ll}
\ns\ds J_1(t)=\dbE_t\int_{t}^{t+\e}\lan  \Big[\sS +\sR(\Th_1+\Th_2)-\big[B^{\top}(\cP_1+\cP_2)+D^{\top}\cP_1(C+D\Th_1+D\Th_2)\big]\Big]X
\\
\ns\ds\qq\qq\qq\qq +\frac 1 2 \sR v +\sR\f -B^{\top}
(\cP_3+\cP_4)-D^{\top}\cP_1(D\f+\si)-D^{\top}\cL_4,v\ran dr.
\ea
\ee
\begin{lemma}\label{Lemma-for-J-1}
Suppose (H1) holds, $X$ solves (\ref{State-equation-Th-1-Th-2}) associated with $(\Th_1,\Th_2,\f)$, and $J_1(t)$ is defined in Lemma \ref{Cost-functional-1}. Then (\ref{Expression-of-J-1}) is true, where $\cP_i$ satisfies (\ref{Riccati-equation-system}).
\end{lemma}

\subsection{A new expression of $J_2$}

In the following, we turn to treating $J_2$. To this end, we introduce
$$\left\{\ba{ll}
\ns\ds dY_0^\e=-\Big[A_\th^{\top}Y_0^\e  + C_\th^{\top} Z_0^\e-F_1^\e\Big]dr+Z_0^\e dW(r),\ \ r\in[t,T],\\
\ns\ds Y_0^\e(T,t)=-G X_0^\e(T)-\wt G  \dbE_tX_0^\e(T),
\ea\right.$$
where $F_1^\e$ is defined in Lemma \ref{Cost-functional-1}. From Proposition 2.1 in \cite{Sun-Yong-2014}, we see that
$$
\ba{ll}
\ns\ds  (Y_0^\e,Z_0^\e) \in  L^2_{\dbF}(\Omega;C([t,T];\dbR^{n}))\times L^2_{\dbF}(t,T;\dbR^{n}),\ \ t\in[0,T).
\ea
$$
Recall $X_0^\e$ in (\ref{Rewrite-of-X-0}), we obtain the following by It\^{o}'s formula,
$$\ba{ll}
\ns\ds d\lan Y_0^\e,X_0^\e \ran=-\lan  A_\th^{\top} Y_0^\e
 + C_\th^{\top} Z_0^\e
  -F_1^\e, X_0^\e \ran dr
  +\lan Z_0^\e, X_0^\e \ran  dW(r)\\
\ns\ds\qq\qq\qq+\lan Y_0^\e, A_\th X_0^\e+B
 vI_{[t,t+\e]} \ran dr+\lan Y_0^\e, C_\th X_0^\e + D vI_{[t,t+\e]}\ran dW(r)\\
\ns\ds\qq \qq\qq+\lan Z_0^\e, C_\th X_0^\e + D vI_{[t,t+\e]}  \ran dr.
\ea
$$
As a result, we then have
\bel{First-one-to-J-1-2}\ba{ll}
\ns\ds \dbE_t\lan -GX_0^\e(T)-\wt G \dbE_tX_0^\e(T),X_0^\e(T)\ran
-\dbE_t\int_t^T \lan F_1^\e,X_0^\e\ran dr\\
\ns\ds =
\dbE_t\int_t^{t+\e}\lan B^{\top} Y_0^\e+ D^{\top} Z_0^\e, v\ran dr.
\ea
\ee
By the decoupling tricks in Lemma \ref{representation-for-Y-Z-*}, we introduce
$$\left\{
\ba{ll}
\ns\ds d\bar\cP_1=-\Big[\bar\cP_1 (A+B\Th_1)+(A+B\Th_1)^{\top}
\bar\cP_1+(C+D\Th_1)^{\top}\bar\cP_1(C+D\Th_1)\\
\ns\ds\qq\qq -\big[Q + S^{\top}\Th_1 +\Th_1^{\top}S+\Th_1^{\top}R\Th_1\big]\Big]ds,\ \ \\
\ns\ds d\bar\cP_2=-\Big\{\bar\cP_2(A+B\Th_1)+(A+B\Th_1)^{\top}\bar\cP_2-\big[\wt Q +\wt S^{\top}\Th_1+\Th_1^{\top}\wt S+\Th_1^{\top}\wt R\Th_1\big]\Big\}ds,\\
\ns\ds d\bar\cP_{3}=-\Big[(A+B\Th_1)^{\top}\bar\cP_3+ \bar\cP_2 BvI_{[t,t+\e]}\Big]ds+\bar\cL_{3}dW(s),\\
d\bar\cP_4=-\Big\{(A+B\Th_1)^{\top}\bar\cP_4+
\big[(C+D\Th_1)^{\top}\bar\cP_1 D +\bar\cP_1 B\big] vI_{[t,t+\e]}\\
\ns\ds\qq\qq +(C+D\Th_1)^{\top}\bar\cL_4\Big\}ds+\bar\cL_4dW(s),\\
\ns\ds \bar\cP_1(T)=-G,\ \bar\cP_2(T)=-\wt G,\ \bar\cP_3(T)=0,\ \bar\cP_4(T)=0.
\ea\right.$$
Moreover, from Lemma \ref{representation-for-Y-Z-*}, the following holds on $[t,T]$,
$$
\ba{ll}
\ns\ds Y_0^\e=\bar\cP_1 X_0^\e +\bar\cP_2 \dbE_tX_0^\e  + \dbE_t\bar\cP_3 +\bar\cP_4,  \ \ Z_0^\e=\bar\cP_1(C+D\Th_1)X_0^\e +\bar\cP_1 DvI_{[t,t+\e]}+\bar\cL_{4}.
%
\ea $$
At this moment, we take a closer look at $(\bar\cP_3,\bar\cL_3),$ $(\bar\cP_4,\bar\cL_4).$ By the uniqueness of BSDEs in Proposition 2.1 of \cite{Sun-Yong-2014}, we have the following equalities
$$\ba{ll}
\ns\ds \bar\cP_3(s)=\wt\cP_3(s)v,\ \ \bar\cL_3(s)=0, \ \ \bar\cP_4(s)=\wt\cP_4(s)v,\ \ \bar\cL_4(s)=0,\ \ s\in[t,T],
\ea
$$
where
$$\left\{\ba{ll}
\ns\ds d\wt\cP_{3}=-\Big[(A+B\Th_1)^{\top}\wt\cP_3+ \bar\cP_2 BI_{[t,t+\e]}\Big]ds, \ \ s\in[t,T],\\
d\wt\cP_4=-\Big\{(A+B\Th_1)^{\top}\wt\cP_4+
\big[(C+D\Th_1)^{\top}\bar\cP_1 D +\bar\cP_1 B\big] I_{[t,t+\e]}\Big\}ds,\ \ s\in[t,T],\\
\ns\ds \wt\cP_3(T)=\wt\cP_4(T)=0.
\ea\right.
$$
Consequently, on $[t,T]$ we conclude that
$$\ba{ll}
\ns\ds B^{\top} Y_0^\e+ D^{\top} Z_0^\e=\big[B^{\top}\bar\cP_1
+D^{\top}\bar\cP_1(C+D\Th_1)\big] X_0^\e+B^{\top}\bar\cP_2\dbE_t X_0^\e\\
\ns\ds\qq \qq\qq\qq +B^{\top} \dbE_t\wt\cP_3+B^{\top}\wt\cP_4+D^{\top}\bar\cP_1D vI_{[t,t+\e]}.
\ea
$$
As a result,
$$\ba{ll}
\ns\ds \dbE_t\int_t^{t+\e}\lan B^{\top} Y_0^\e+ D^{\top} Z_0^\e, v\ran dr\\
\ns\ds
=\dbE_t\int_t^{t+\e}\lan B^{\top}\big[\bar\cP_1+\bar\cP_2+D^{\top}\bar\cP_1(C+D\Th_1)\big]X_0^\e+
B^{\top}[\wt\cP_3+\wt\cP_4]+D^{\top}\bar\cP_1Dv,v\ran dr.
\ea
$$
By the estimate of $X_0^\e$, for almost $t\in[0,T),$
$$\ba{ll}
\ns\ds  \dbE_t\int_t^{t+\e}\lan B^{\top}\big[\bar\cP_1+\bar\cP_2+D^{\top}\bar\cP_1(C+D\Th_1)\big]X_0^\e,v\ran dr=o(\e).
\ea
$$
From the equations of $(\wt\cP_3,\wt\cP_4)$,
$$\ba{ll}
\ns\ds \sup_{t\in[t,t+\e]}\big[|\wt\cP_3(t)|^2+|\wt\cP_4(t)|^2\big]=o(\e).
\ea
$$
To sum up, by the definition of $J_2$ and (\ref{First-one-to-J-1-2}), for almost $t\in[0,T)$ we deduce that
\bel{Expression-of-J-2}\ba{ll}
\ns\ds  J_2(t)=\frac \e 2 \lan D(t)^{\top}\bar\cP_1(t) D(t)v,v\ran+o(\e).
\ea
\ee
\begin{lemma}\label{Cost-functional-2}
Suppose (H1) holds, $X_0^\e$ is in (\ref{Rewrite-of-X-0}) associated with $(\Th_1,\Th_2,\f)$, and $J_2(t)$ is defined in Lemma \ref{Cost-functional-1}. Then (\ref{Expression-of-J-2}) is true.
\end{lemma}

\subsection{Proofs of the main results}\label{Main-proof-subsection}

We are in the position to give the proofs of the main results in Section 3.

To begin with, we give the proof of Theorem \ref{Equivalence-open-1}.

\begin{proof}
In Lemma \ref{Cost-functional-1}, Lemma \ref{Lemma-for-J-1}, Lemma \ref{Cost-functional-2}, we take $\Th_1\equiv\Th_2\equiv0$. Hence for the notations in (\ref{Definition-of-u-any}), $u\equiv\f$ and
$$\left\{\ba{ll}
\ns\ds J_1(t)=\dbE_t\int_{t}^{t+\e}\lan  \Big[\sS-\big[B^{\top}(P_1+P_2)+D^{\top}P_1C\big]\Big]X
+\frac 1 2 \sR v\\
\ns\ds\qq\qq\q  +\sR u -B^{\top}
(P_3+P_4)-D^{\top}P_1(D u+\si)-D^{\top}L_4,v\ran dr,\\
\ns\ds J_2(t)=\frac \e 2\lan D(t)^{\top}P_1(t)D(t)v,v\ran +o(\e),
\ea\right.
$$
where $P_i$, $i=1,2$, $(P_j,L_j)$, $j=3,4,$ satisfies (\ref{Systems-for-open-only}). Moreover, for any $t\in[0,T),$ by Remark \ref{estimates},
$$\ba{ll}
\ns\ds \dbE_t\int_t^{t+\e}\lan \sS^{\top}v,X_0^\e\ran ds=o(\e).
\ea
$$
We set out to define $\bar X$ the state process associated with $\bar u$,  $u^{v,\e}:=\bar u+vI_{[t,t+\e]}$, and for any $t\in[0,T)$
\bel{Notation-proof-main-result-1}\left\{\!\! \ba{ll}
\ns\ds \sD_0(t):=\lim_{\e\rightarrow0}\frac{1}{2\e}  \int_t^{t+\e}
\big[\sR(s)-D (s)^{\top}P_1(s)D(s)\big]ds, \\
\ns\ds \sH_0(t):=\lim_{\e\rightarrow0}\frac 1 \e
\dbE_t\int_t^{t+\e}\Big[\sS(s) \bar X(s)+\sR(s)\bar u(s)-B(s)^{\top}\bar M(s,s)-D(s)^{\top}\bar N(s)\Big]ds
\ea\right.\ee
with $(\bar M,\bar N)$ in (\ref{Arbitrary-s-Y-Z-general-2}) corresponding to $\bar u.$
To sum up, $u\equiv\bar u=\bar\f$ is an equilibrium control associated with $x_0$ if and only if for any $t\in[0,T),$ $v\in L^2_{\cF_t}(\Omega;\dbR^m)$,
$$\ba{ll}
\ns\ds 0\leq \lim_{\e\rightarrow0}\frac{J(t,\bar X(t);u^{v,\e}(\cd))-J\big(t,\bar X(t);\bar u(\cd)\big)}{\e}
=\lan \sD_0(t) v,v\ran+\lan \sH_0(t),v\ran.
\ea$$
Given $t\in[0,T)$, this holds if and only if both $\sH_0(t)=0$ and $\sD_0(t)\geq0$. Since both $\sR$ and $P_1$ are bounded and deterministic, we thus know that
$$
0\leq \sR(t)-D(t)^{\top}P_1(t)D(t),\ \ t\in[0,T].\ \ a.e.
$$
If $\sH_0(t)=0$, then by Lemma 3.4 in \cite{Hu-Jin-Zhou-2017}, above (\ref{necessary-open-loop-2}) holds. Conversely, if (\ref{necessary-open-loop-2}) is true, we immediately obtain $\sH_0(t)=0$.
\end{proof}

\medskip

Next we present the proof of Theorem \ref{Th-closed-loop-representation}.

\begin{proof}
In Lemma \ref{Cost-functional-1}, Lemma \ref{Lemma-for-J-1}, Lemma \ref{Cost-functional-2}, we take $\Th_1\equiv0$. Hence for the notations in (\ref{Definition-of-u-any}), we have $u\equiv\Th_2 X+\f$ and
$$\left\{\ba{ll}
\ns\ds J_1(t)=\dbE_t\int_{t}^{t+\e}\lan  \Big[\sS +\sR\Th_2-\big[B^{\top}(\cP_1+\cP_2)+D^{\top}\cP_1(C+D\Th_2)\big]\Big]X
+\frac 1 2 \sR v\\
\ns\ds\qq\qq\q  +\sR\f -B^{\top}
(\cP_3+\cP_4)-D^{\top}\cP_1(D\f+\si)-D^{\top}\cL_4,v\ran dr,\\
\ns\ds J_2(t)=\frac \e 2\lan D(t)^{\top}P_1(t)D(t)v,v\ran +o(\e),
\ea\right.
$$
where $\cP_i$, $i=1,2$, $(\cP_j,\cL_j)$, $j=3,4,$ satisfies (\ref{System-closed-loop-representation}). Moreover, by Remark \ref{estimates},
$$\ba{ll}
\ns\ds \dbE_t\int_t^{t+\e}\lan \sS^{\top}v,X_0^\e\ran ds=o(\e),\ \ t\in[0,T).
\ea
$$
For open-loop equilibrium strategy pair $(\Th^*,\f^*)$ and associated equilibrium control $u^*$, we define $X^*$ the corresponding state process as,
$$\left\{\ba{ll}
\ns\ds dX^*=\big[(A+B\Th^*) X^*+B\f^*+b \big]ds +\big[(C+D\Th^*) X^*+D\f^*+\si\big]dW(s),\\
\ns\ds X^*(0)=x_0,
\ea\right.$$
and perturbed control $u^{v,\e}:=\Th^* X^*+\f^*+vI_{[t,t+\e]}$. Moreover, for $(\cM^*,\cN^*)$ in (\ref{Definition-c-M-N}) corresponding to $u^*,$ let
$$\ba{ll}
\ns\ds \sH_1(t):=\lim_{\e\rightarrow0}\frac 1 \e
\dbE_t\int_t^{t+\e}\Big[\sS(s) X^*(s)+\sR(s) u^*(s)-B^{\top}\cM^*(s,s)-D^{\top}\cN^*(s)\Big]ds.
\ea$$
To sum up, $u^*=\Th^*X^*+\f^*$ is an equilibrium control associated with $x_0\in\dbR^n$ if and only if for any $t\in[0,T],$ $v\in L^2_{\cF_t}(\Omega;\dbR^m)$,
\bel{Inequality-condition-for-J-0}\ba{ll}
\ns\ds 0\leq \lan \sD_0(t) v,v\ran+\lan \sH_1(t),v\ran,
\ea\ee
where $\sD_0$ is in (\ref{Notation-proof-main-result-1}). Given $t\in[0,T)$, this holds if and only if both $\sH_1(t)=0$ and $\sD_0(t)\geq0$. Since both $\sR$ and $P_1$ are bounded and deterministic,
$$
0\leq  \sR(t)-D(t)^{\top}P_1(t)D(t),\ \ t\in[0,T].\ \ a.e.
$$
$\Longrightarrow$  If $\sH_1(t)=0$, then by Lemma 3.4 in \cite{Hu-Jin-Zhou-2017}, for almost $s\in[0,T],$ we have
\bel{Equivalence-2}\ba{ll}
\ns\ds
0=\sS  X^* +\sR u^*-B^{\top}\cM^* -D^{\top}\cN^* \\
\ns\ds\q =\Big[\sS +\sR\Th^*-\big[B^{\top}(\cP_1^*+\cP_2^*)+D^{\top}\cP_1^*(C+D\Th^*)\big]\Big]X^*
 \\
\ns\ds\qq\q +\sR\f^* -B^{\top}
(\cP_3^*+\cP_4^*)-D^{\top}\cP_1^*(D\f^*+\si)-D^{\top}\cL_4^*.
\ea
\ee
Notice that (\ref{Equivalence-2}) holds for any $x_0\in\dbR^n$. We choose $x_0=0$, and denote the state process by $X^*_0$. As a result,
$$\ba{ll}
\ns\ds \Big[\big[\sR- D^{\top}\cP_1^* D\big]\Th^*-
 B^{\top}\big[\cP_1^*+\cP_2^*\big]-D^{\top}\cP_1^* C+\sS\Big](X^*-X_0^*)=0.
\ea
$$
At this moment, given $I\in\dbR^{n\times n}$ the unit matrix, we consider the following equation
\bel{Standard-unit-SDE}\left\{ \ba{ll}
\ns\ds d\sX=(A+B\Th^*) \sX ds +(C+D\Th^*)\sX  dW(s),\ \ s\in[0,T],\\
\ns\ds \sX(0)=I,
\ea\right. \ee
the solvability of which is easy to see. Moreover, $ \sX ^{-1}$ also exists. By the standard theory of SDEs,
$$\dbP\big\{\omega\in\Omega;\ \sX(t,\omega)x=X^*(t,\omega)-X_0^*(t,\omega),\ \forall t\in[0,T]\big\}=1.$$
Using the existence of $\sX^{-1}$, it is easy to see above (\ref{One-equality-for-Theta}).

$\Longleftarrow$ In this case, it is easy to see (\ref{Equivalence-2}) with $u^*:=\Th^*X^*+\f^*$. Consequently, the conclusion is followed by (\ref{Inequality-condition-for-J-0}), (\ref{necessity-open-loop-1}) and the fact of $\sH_1(t)=0$.
\end{proof}

\medskip

At last, we show the proof of Theorem \ref{Th-closed-loop-strategy}.
\begin{proof}
In Lemma \ref{Cost-functional-1}, Lemma \ref{Lemma-for-J-1}, Lemma \ref{Cost-functional-2}, we take $\Th_2\equiv0$. Hence for the notations in (\ref{Definition-of-u-any}), $u\equiv\Th_1 X+\f$ and
$$
\left\{\ba{ll}
\ns\ds J_1(t)=\dbE_t\int_{t}^{t+\e}\lan  \Big[\sS +\sR\Th_1 -\big[B^{\top}(\sP_1+\sP_2)+D^{\top}\sP_1(C+D\Th_1)\big]\Big]X
+\frac 1 2 \sR v\\
\ns\ds\qq\qq\q  +\sR\f -B^{\top}
(\sP_3+\sP_4)-D^{\top}\sP_1(D\f+\si)-D^{\top}\sL_4,v\ran dr,\\
\ns\ds J_2(t)=\frac \e 2\lan D(t)^{\top}\sP_1(t)D(t)v,v\ran +o(\e),
\ea\right.
$$
where $\sP_i$, $i=1,2$, $(\sP_j,\sL_j)$, $j=3,4,$ satisfies (\ref{System-of-closed-loop-strategy}). Moreover, in view of Remark \ref{estimates}, it is straightforward to get
$$\ba{ll}
\ns\ds \dbE_t\int_t^{t+\e}\lan (\sS^{\top}+\Th_1^{\top}\sR)v,X_0^\e\ran ds=o(\e),\ \ t\in[0,T).
\ea
$$
For closed-loop equilibrium strategy pair $(\Th^*,\f^*)$ in the sense of Definition \ref{Definition-3} and associated equilibrium control $u^*:=\Th^* X^*+\f^*$, we define $X^*$ the corresponding state process as,
$$\left\{\ba{ll}
\ns\ds dX^*=\big[(A+B\Th^*) X^*+B\f^*+b \big]ds +\big[(C+D\Th^*) X^*+D\f^*+\si\big]dW(s),\\
\ns\ds X^*(0)=x_0,
\ea\right.$$
and perturbed control variable $u^{v,\e}:=\Th^* X^{v,\e}+\f^*+vI_{[t,t+\e]}$. In addition, for $(\sM^*,\sN^*)$ in (\ref{Definitions-of-s-M-s-N}) corresponding to $u^*,$ we denote by
$$\left\{\ba{ll}
\ns\ds \sH_2(t):=\lim_{\e\rightarrow0}\frac 1 \e
\dbE_t\int_t^{t+\e}\Big[\sS(s) X^*(s)+\sR(s) u^*(s)-B^{\top}\sM^*(s,s)-D^{\top}\sN^*(s)\Big]ds,\\
\ns\ds \sD_1(t):=\lim_{\e\rightarrow0}\frac{1}{2\e}\int_t^{t+\e}\big[\sR(s)
-D(s)^{\top}\sP_1^*(s)D(s)\big]ds.
\ea\right.$$
To sum up, $u^*=\Th^*X^*+\f^*$ is a closed-loop equilibrium control associated with $x_0\in\dbR^n$ if and only if for any $t\in[0,T],$ $v\in L^2_{\cF_t}(\Omega;\dbR^m)$,
\bel{Inequality-condition-for-J-1}\ba{ll}
\ns\ds 0\leq \lan \sD_1(t) v,v\ran+\lan \sH_2(t),v\ran.
\ea\ee
Given $t\in[0,T)$, this holds if and only if both $\sH_2(t)=0$ and $\sD_1(t)\geq0$.

$\Longrightarrow$  Given equilibrium strategy pair $(\Th^*,\f^*)$, we conclude that
$\sP_1^*$ is bounded and deterministic. Recall the requirement on $\sR$, it is clear that
\bel{Second-order-equilibrium-closed-loop-strategy}\ba{ll}
\ns\ds
0\leq \sR(t)-D(t)^{\top}\sP_1^*(t)D(t),\ \ t\in[0,T].\ \ a.e.
\ea\ee
If $\sH_2(t)=0$, then by Lemma 3.4 in \cite{Hu-Jin-Zhou-2017}, for almost $s\in[0,T],$ we have
\bel{Equivalence-3}\ba{ll}
\ns\ds
0=\sS  X^* +\sR u^*-B^{\top}\sM^* -D^{\top}\sN^* \\
\ns\ds\q =\Big[\sS +\sR\Th^*-\big[B^{\top}(\sP_1^*+\sP_2^*)+D^{\top}\sP_1^*(C+D\Th^*)\big]\Big]X^*
 \\
\ns\ds\qq\q +\sR\f^* -B^{\top}
(\sP_3^*+\sP_4^*)-D^{\top}\sP_1^*(D\f^*+\si)-D^{\top}\sL_4^*.
\ea
\ee
Notice that (\ref{Equivalence-3}) holds for any $x_0\in\dbR^n$. We choose $x_0=0$, and denote the state process by $X^*_0$. As a result,
$$\ba{ll}
\ns\ds \Big[\big[\sR- D^{\top}\sP_1^* D\big]\Th^*-
 B^{\top}\big[\sP_1^*+\sP_2^*\big]-D^{\top}\sP_1^* C+\sS\Big](X^*-X_0^*)=0.
\ea
$$
As in Theorem \ref{Th-closed-loop-representation}, we introduce $\sX$ satisfying (\ref{Standard-unit-SDE}), and therefore obtain (\ref{Conclusion-of-closed-strategy}) by following the same spirit of that in Theorem \ref{Th-closed-loop-representation}.

$\Longleftarrow$ In this case, it is easy to see (\ref{Equivalence-2}) with $u^*:=\Th^*X^*+\f^*$. Consequently, the conclusion is followed by (\ref{Inequality-condition-for-J-1}), (\ref{Second-order-equilibrium-closed-loop-strategy}) and the fact of $\sH_1(t)=0$.
\end{proof}

\section{Concluding remarks}

In the Markovian setting, a unified approach by variational idea is developed to build the characterizations for three notions, i.e., closed-loop equilibrium controls/strategies, open-loop equilibrium controls, as well as the closed-loop representations of open-loop equilibrium controls. The intrinsic differences among different equilibrium controls are also revealed clearly and deeply. Related studies with random coefficients or in mean-field setting are under consideration. We hope to do some relevant research in future.

\end{document}